\documentclass[12 pt]{amsart}
\usepackage[latin1]{inputenc}
\usepackage[cyr]{aeguill}
\usepackage[english]{babel}
\usepackage{amsxtra}
\usepackage[dvips]{graphicx}
\usepackage{amsmath,amsthm, amssymb,amsopn,pst-tree,stmaryrd,wasysym}
\usepackage[all]{xy}
\newtheorem{theoreme}{\sc Theorem}
\newtheorem{proposition}{\sc Proposition}
\newtheorem*{lemme}{\sc Lemma}
\newtheorem*{corollaire}{\sc Corollary}
\theoremstyle{remark}
\newtheorem{remarque}{\it Remark}
\newtheorem{notation}{\it Notation}

\newtheorem{definition}{\sc Definition}
\font\tmsb=msbm10 at12pt
\font\smsb=msbm7
\font\ssmsb=msbm5
\newfam\msbfam
\textfont\msbfam=\tmsb
\scriptfont\msbfam=\smsb
\scriptscriptfont\msbfam=\ssmsb
\def \RM{\mathbb {R}}%        corps des reels
\def \NM{\mathbb{N}}%        entiers naturels
\def \ZM{\mathbb{Z}}%        entiers relatifs
\def \CM{\mathbb{C}}%        nombres complexes
\def \QM{\mathbb{Q}}%        nombres rationnels

\def\2n{(\CM^{2n},0)}

\newcommand{\mat}[4]{\left(\!\begin{array}{cc} #1 & #2 \\ #3 & #4 \end{array}
\!\right)}
\newcommand{\app}[5]{\begin{array}{ccccc} #1 : & #2  & \longrightarrow & #3 
\\ & #4 & \mapsto & #5 \end{array}}
\medskip
\title[The ring of quasimodular forms for a cocompact group]
{\sl The ring of quasimodular forms for a \\  cocompact  group}
\author{\sl Najib Ouled Azaiez }
\date{Mars 2006}
\address{Institut Math{\'e}matiques de Jussieu, 175 Rue de
Chevaleret,75013 Paris}
\email{ouled@math.jussieu.fr}
\begin{document}
\maketitle
%%%%%%%%%%%%%%%%%%%%%%%%%%%%%%%%%%%%%%%%%%%%%%%%%%%%%%%%
%%%%%%%%%%%%%%%%%%%%%%%%%%%%%%%%%%%%%%%%%%%%%%%%%%%%%%%%%
%%%%%%%%%%%%%%%%%% ABSTRACT%%%%%%%%%%%%%%%%%%%%%%%%%%%%%%
%%%%%%%%%%%%%%%%%%%%%%%%%%%%%%%%%%%%%%%%%%%%%%%%%%%%%%%%%
\begin{abstract}
We describe the  additive  structure of the
graded ring $\widetilde{M}_*$ of 
quasimodular forms over any discrete and 
cocompact group $\Gamma \subset \rm{PSL}(2, \RM).$
We show that  this ring is never finitely generated.
We calculate the exact number of new generators
in each weight $k$. This number is constant
for $k$ sufficiently large and equals $\dim_{\CM}(I / 
I \cap \widetilde{I}^2),$ where $I$ and $\widetilde{I}$
are the ideals of modular forms and  quasimodular forms,
respectively, of 
positive weight. We show that $\widetilde{M}_*$
is contained in some finitely generated ring
$\widetilde{R}_*$ of meromorphic quasimodular forms
with $\dim \widetilde{R}_k = O(k^2),$ i.e. the same
order of growth as $\widetilde{M}_*.$ 
\end{abstract}
\section{Introduction}
Kaneko and Zagier introduced the notion of
quasimodular forms in \cite{Kan-Zag}.  The structure 
of $\widetilde{M}_*(\Gamma_1)$ (where $\Gamma_1
= \rm{PSL}(2, \ZM)$  is the classical modular group) 
was given in \cite{Kan-Zag}, in which it is proved that 
$\widetilde{M}_*(\Gamma_1) = \mathbb{C}[E_2,E_4,E_6],$ with
$E_2,E_4$ and $E_6$ being the Eisenstein 
series of weights $2,4$ and $6$ respectively. 

We study the ring of quasimodular forms over
discrete and cocompact subgroups of $\rm{PSL}
(2,\RM).$ In the second and third section, we 
derive some general properties of quasimodular 
forms over discrete and cocompact subgroups of 
$\rm{PSL}(2, \RM),$ following $[14],$ $[4]$ and $[15].$
In the end of the third section, we give
an additive structure theorem of rings of quasimodular
forms and a $\rm{sl}_2(\CM)$-module structure theorem
for the ring of quasimodular forms.
In the fourth section, we give a cocompact/no cocompact dichotomy
(see Theorem $4$) which characterizes cocompact modular groups
in terms of their spaces of quasimodular forms of weight $2.$
In the  fifth section, we describe
our principal results  which are
Theorem $5$ and its corollary. In Theorem $6,$ we
describe the additive structure of the differential 
closure of any ring
$\mathcal{M}$ generated by holomorphic or
meromorphic modular forms of positive weights
over any discrete and finite covolume subgroup
of $\rm{PSL}(2, \RM).$ We prove also that the
differential closure of $\mathcal{M}$ is
never finitely generated.  In the last section, we
prove the existence of quasimodular forms of 
weight $2$ with prescribed poles. We use this result
to construct finitely generated rings of meromorphic
quasimodular forms with positive weights, over 
cocompact groups (see Theorem $10$ and its corollary).
In the sixth section, we give an algebraic characterization
of cocompact  groups, in terms of their modular
forms rings (see Theorem $8$).
\section{General properties of quasimodular forms}
In this section, we recall definitions and general
properties of quasimodular forms, given by Kaneko
and Zagier in $[4].$ We give corollaries 
of sevral results in $[4]$ and new proofs for other 
results in this paper.

We consider a discrete and finite covolume subgroup
$\Gamma$ of $\rm{PSL}(2,\RM).$ We give the definition
of modular forms, quasimodular forms, almost holomorphic
modular forms and modular stacks, over the group $\Gamma.$
We denote  by  $\mathcal{H}$ the upper half plane and by
$y$ the imaginary part of $z \in \mathcal{H}.$

\begin{definition}
A modular form of weight $k$ over $\Gamma$ is
a holomorphic map $f$ in  $\mathcal{H}$ with
moderate growth, such that : 
\begin{equation}\label{def:mod} 
(cz +d)^{-k} f(\frac{a z + b} {c z + d}) = f(z), \quad\forall
\mat{a} {b} {c} {d}  \in \Gamma \mbox{ } \mbox{and} \mbox{ }
z \in \mathcal{H}.\end{equation} 
\end{definition}
\begin{definition}
A quasimodular form $f$ of weight $k$ and depth $\leq p$
over $\Gamma,$ is a holomorphic function $f$ in 
$\mathcal{H}$ with moderate growth, such that for any
$z \in \mathcal{H},$ the map :
$$\begin{array}{cccc}
\Gamma & \longrightarrow & \CM \\
\mat{a} {b} {c} {d} & \mapsto & (c z + d)^{-k} 
f(\frac{a z + b} {c z + d}),\end{array}$$ is a 
polynomial of degree $\leq p$ in $\frac{c} {c z + d}$ 
with functions defined on $\mathcal{H}$ as coefficients.
We can write : 
\begin{equation}\label{def:qua}
(c z + d)^{-k} f(\frac{a z + b} {c z +  d})
= \sum_{j = 0}^{p} f_j(z) (\frac{c} {c z + d})^j,\quad \forall
z \in \mathcal{H},\end{equation}
with map $f_j : \mathcal{H} \longrightarrow
\CM \quad (j=0,\cdots,p).$
\end{definition}
\begin{remarque}
This definition which is different from the one given in \cite{Kan-Zag} 
was proposed by Werner Nahm. The equivalence between this definition and 
the one given in \cite{Kan-Zag} is a consequence of Theorem $1.$
\end{remarque}
\begin{definition}\label{def:pres}
An almost holomorphic modular form $F$ of weight $k$ and
depth $\leq p$ over $\Gamma$ is a polynomial in
$\frac{1} {y}$ of degree $\leq p$ whose coefficients
are holomorphic maps on $\mathcal{H}$ with
moderate growth, such that  (\ref{def:mod}) holds for any
$\mat{a} {b} {c} {d} \in \Gamma$ and $z \in \mathcal{H}.$
We can write :
$$F(z) = f_0(z) + \frac{f_1(z)} {z - \overline{z}} +
\cdots + \frac{f_p(z)} {(z - \overline{z})^p},$$ 
with holomorphic maps $(f_i),$ because $y = \frac{z -
\overline{z}} {2 i}.$
\end{definition}
This way of writing $F$ as polynomial in $\frac{1} {z - 
\overline{z}}$ is more useful for making the next calculations.
\begin{definition}
A modular stack of weight $k$ and depth $\leq p$ is
a holomorphic map :
$$\app{E} {\mathcal{H}} {\bigoplus_{l = 0}^{\infty} \CM} 
{z} {(f_0(z),f_1(z),\cdots)}$$
with moderate growth such that the maps $f_l$ satisfy 
$f_l = 0$ for $l > p$ and the functional equation :
\begin{equation}\label{mod:field}
(c z + d)^{-k + 2 l} f_l( \frac{a z + b} {c z + d})
 = \sum_{j \geq l} \binom{j} {l} f_j(z) (\frac{c} 
{c z + d})^{j - l}.\end{equation}
\end{definition}
\begin{notation}
We denote by $M_* = \bigoplus_{k \geq 0} M_k$ (respectively
$\widetilde{M}_* = \bigoplus_{k \geq 0} \widetilde{M}_k,$
$\widehat{M}_* = \bigoplus_{k \geq 0} \widehat{M}_k,
\overrightarrow{M}_* = \bigoplus_{k \geq 0} \overrightarrow{M}_k$)
the graded rings of modular forms (respectively quasimodular
forms, almost holomorphic modular forms and modular stacks).
We denote by $\widetilde{M}_*^{(\leq p)},\widehat{M}_*^{(\leq p)},
\overrightarrow{M}_*^{(\leq p)}$
the subspaces of quasimodular forms (respectively almost 
holomorphic modular forms and modular fields) of depth 
$\leq p$ over a  certain group $\Gamma.$
\end{notation}
\begin{theoreme}\label{Kan-Z}
Let $\Gamma \subset \rm{PSL}(2 , \RM)$ be a discrete
subgroup of finite covolume and $p$ a positive integer. 
We have the isomorphisms :
$$\begin{array}{llcll}
\widetilde{M}_*^{(\leq p)} &  \simeq &  
\overrightarrow{M}_*^{(\leq p)}  & \simeq &
\widehat{M}_*^{(\leq p)} \\   f & \mapsto &(f_0,\cdots,f_p) 
& \mapsto &\sum_{j=0}^{p} \frac{f_j(z)} {(z - \overline{z})^j}, 
\end{array}$$ where the sequence of coefficients $(f_j)$ are 
associated to $f$ according to $(2).$ The inverse map
is given by $(f_0,\cdots,f_p) \longrightarrow f_0.$
\end{theoreme}
\begin{remarque}
This theorem implies that an almost holomorphic modular 
form and a modular stack  are determined by  their first
coefficient $f_0,$ or  first coordinate $f_0$
respectively.
\end{remarque}
\section{The additive and $\rm{sl}_2(\CM)$-modules
structure of rings of  quasimodular forms }
There  exists three derivation operators on the spaces
of quasimodular forms.  By the isomorphisms of Theorem $1,$
we get the corresponding operators on the other spaces.
We check that there exists a representation of the Lie algebra
$\rm{sl}_2(\CM)$ on the spaces $\widetilde{M}_*,$
$\widehat{M}_*$ and $\overrightarrow{M}_*$ of quasimodular forms,
almost holomorphic modular forms and modular stacks.
\begin{proposition}
The operator $D$ of derivation, with respect to $z$ acts on
the space of quasimodular forms. This operator increases
the weight by $2$ and the depth by $1.$ For any
$k \geq 0$ and $p \geq 0$ we have:
$$D: \widetilde{M}_k^{(\leq p)} \longrightarrow 
\widetilde{M}_{k+2}^{(\leq p + 1)}.$$
\end{proposition}
\begin{proof}
Let $f \in \widetilde{M}_k ^{(\leq p)}.$ By definition we
have: $$(c \; z + d)^{-k}  
f (\frac{ a \; z + b} {c \; z + d}) = \sum_{ 0 \leq j \leq p} f_j(z) \; 
(\frac{c} {c \; z + d})^j,$$ with holomorphic maps  $f_j$.   
So:
\[
\begin{array}{rcl} &  & (c \; z + d)^{-k-2} \;  f'( \frac{a \; z + b} 
{c \; z + d} )  \\ & & \\ = &  
&  D[(c \; z +d)^{-k} \; f(\frac{ a \; z +b} {c \; z +d })]+ k c \; 
( c\; z + d)^{-k-1} \; f( \frac{ a 
\; z + b} { c \; z + d} ) \\ &  & \\= &  &  D[\sum_{ 0 \leq j \leq p} 
f_j(z) \; ( \frac{c} 
{ c \; z + d} )^j] + \frac{k c} {c  \; z + d}\sum_{ 0 \leq j \leq p}  
f_j(z)  \; (\frac{c} 
{ c \; z + d})^j\\ & &  \\ = &  &  \sum_{0 \leq j \leq p+1} [f'_j(z) + 
(k-j+1) f_{j-1}(z)] \; 
(\frac{c} {c  \; z + d})^j.  \end{array} \]
(with $f_{-1} \equiv f_{p+1} \equiv 0 $ ). So the weight of $f'$ 
is $k+2$ and its depth is $\leq (p+1).$ 
\end{proof}

\begin{proposition}\label{pro:del:pres}
If  $f \in \widetilde{M}_k^{(\leq p)}$ is a quasimodular
form and  $F(z)=f_0(z) + \frac{f_1(z)} {z - 
\overline{z}} + \cdots + \frac{f_p(z)} {(z - \overline{z})^p}$
with $f_0 = f$ is the almost holomorphic modular form
which corresponds to $f$, then any $f_l$ is a quasimodular
form of weight $k - 2l$ and depth $\leq p - l.$ In particular,
we have a map $\delta : \widetilde{M}_k \longrightarrow 
\widetilde{M}_{k-2}$ which maps  $\widetilde{M}_k ^{( \leq p)}$ 
to $\widetilde{M}_{k-2}^{(\leq p-1)}$ for any $p.$ This map is 
given by $f=f_0 \mapsto f_1$ and it has the following
properties:

(i) The kernel of the map $\delta:\widetilde{M}_k \longrightarrow 
\widetilde{M}_{k-2}$ is the space $M_k$.

(ii) If $f(z)$ is a quasimodular form, then the almost
holomorphic modular form  $F$ associated to $f$ is given 
 by: $F(z) = \sum_{n=0}^{\infty} \frac{(\delta^n f)(z)} 
{n!~(z- \overline{z})^{n}}.$
\end{proposition}
\begin{remarque}{\rm  The sum in (ii) is finite 
because if the depth of $f$ is $\leq p$ then
$\delta^n (f)=0$ for  $n >p$ . In fact, since
$\widetilde{M}_k^{(\leq 0)} = M_k$ vanishes for 
$k < 0$, we see that the depth of  
a quasimodular form $f$ of weight $k$ is at most
equal to $\frac{k} {2}.$}   
\end{remarque}
\begin{proof}
The second part $(ii)$ is clear by  using Definition 
$\ref{def:pres}$ and Theorem $\ref{Kan-Z}$. The first part $(i)$ 
is a consequence of this since :
$\delta (f) =0 \Leftrightarrow \delta ^n (f) =0 \;(\forall n \geq 1)  
\Leftrightarrow f = F \Leftrightarrow F$ holomorphic.
\end{proof}
\begin{corollaire}
Let $k \geq 0$, $f \in \widetilde{M}_k ^{(\leq p)}$ and $F = f_0 + 
\frac{f_1} {z - \overline{z}} + \cdots + \frac{f_p} {(z - 
\overline{z})^p} \vspace{1mm}$ the associated almost holomorphic 
modular form.Then, we have $f_p \in M_{k - 2 p}$ and more generally 
$f_j \in \widetilde{M}_{k - 2 \; j} ^{(\leq p - j)}.$
\end{corollaire}
\begin{proof}
By properties of  $\delta,$ it is clear that $f_j \in 
\widetilde{M}_{k-2 \; j} ^{(\leq p-j)}.$ In
particular $f_p \in \widetilde{M}_{k- 2 \; p} ^{(\leq 0)}$. 
Since a quasimodular form of depth $0$ is modular, we deduce
that $f_p \in M_{k - 2 \; p}.$
\end{proof}
\begin{definition}
Let  $ H :\widetilde{M}_* \longrightarrow \widetilde{M}_*$ be 
the operator which associates $k f$ to any 
quasimodular form $f$ of weight $k,$ 
i.e. $H(f) = k \; f$.
\end{definition}

\begin{proposition}
The operators $D$, $\delta$ and  $H$ satisfy the  relations :  
$$ \begin{array}{cc} i) & [H,D] = 2 
\; D. \\ ii)& [H,\delta] = - 2 \; \delta. \\ iii)& [\delta , D] = H. 
\end{array}$$ 
In other words, we have a representation of the Lie algebra
$\rm{sl(2,\mathbb{C})}$ over the spaces
$\widetilde{M}_*$, $\widehat{M}_*$ and $\overrightarrow{M}_*.$
\end{proposition}
\begin{proof}
The points $(i)$ and  $(ii)$ follow from the fact 
applying  $D$ or $\delta,$ increases or decreases
the weight by $2.$ 

To prove  $(iii),$ we compute the braket $[\delta,D]$ over
spaces of modular stacks.
By using Theorem $\ref{Kan-Z},$ we obtain the corresponding result
over spaces of quasimodular forms. 
It is easy to check 
that, for a modular stack  of weight $k$
and depth $\leq p,$ we have the property:    
$$\delta (f_0, \cdots ,f_p) = (f_1,2 f_2 , \cdots , j f_j,\cdots ).$$
We deduce from  Proposition $2$   that: 
$$\begin{array}{llll}  D \delta (f_0, \cdots, f_j, \cdots ,f_p)   &  =   
D(f_1,2 f_2 , \cdots  ,p f_p)\\  &  = (f'_1,\cdots,(j+1) f'_{j+1} + 
(k -1 -j)j f_j, \cdots ).\end{array}$$ 
On the other hand, $$D(f_0,\cdots,f_j,\cdots,f_p) =
(f'_0,f'_1 + k \; f_0,\cdots,f'_j + (k-j+1)f_{j-1}, \cdots ).$$ 
So, 
$$ \delta D(f_0,\cdots ,f_p)= (f'_1+ k \; f_0,\cdots,(j+1) f'_{j+1} + 
(j+1)(k-j)f_j,\cdots ).$$ By taking the difference of the 
two preceding equations (giving $\delta D$ and  $D \delta$ for a modular stack),  
we find:
$$ [\delta ,D] (f_0 , \cdots ,f_p) = 
( k \; f_0 , \cdots,k f_j ,\cdots ) = k(f_0,\cdots,f_p),$$     
which implies (by isomorphisms of Theorem $\ref{Kan-Z}$) 
the property:  $$[\delta , D](f)=  H(f).$$ 
\end{proof}

Let $\mathcal{U}$ be the  universal envelopping algebra of 
$\rm{sl}_2(\CM)$.
We compute next  the class of the    
operator $\delta^n D^n$ modulo $\mathcal{U} \delta,$
for any $n \in \mathbb{N}.$ 
\begin{proposition}\label{pro:cal:ope}
The class of the operator $\delta ^n D^n$ modulo $\mathcal{U} 
\delta$  is given by :
$$(\delta ^n \; D^n)  \equiv 
n! \; \prod_{j=0}^{n-1} (H+j) (\mbox{mod} \mbox{ }
\mathcal{U} \delta) \quad\quad.$$
\end{proposition}
\begin{proof}
By Proposition $3$, we know that 
$\delta  D \equiv  H  (\rm{mod} \; \mathcal{U} \delta).$
Let $j > 1,$ and  suppose that $\delta ^{j-1} 
\; D^{j-1}  \equiv  P_{j-1}(H) (\rm{mod} \; \mathcal{U} \delta)$  
where
$P_{j-1}$ is a polynomial of  degree $(j-1).$ We have:
$$\begin{array}{llcll} \delta^j  D^j  & \equiv  \delta^{j-1}  (\delta D) \; 
D^{j-1} \equiv  \delta^{j-1}
(D\delta + H)D^{j-1} \vspace{1mm}\\ &  \equiv  \delta^{j-2}(\delta D) \delta 
D^{j-1} + \delta^{j-1}  H D^{j-1}
\vspace{1mm}  \\  &  \equiv   \delta ^{j-2} \; D \delta^2  D^{j-1} + 
\delta^{j-2}  H \delta  D^{j-1} + 
\delta^{j-1} H  D^{j-1} \vspace{1mm}\\ &  =  \cdots \\
 & \equiv \delta D \;\delta^{j-1} \; D^{j-1} + \sum_{n=1}^{j-1} \delta^n \; H 
\; \delta^{j-1-n} \; D^{j-1} (\rm{mod} \; \mathcal{U} \delta)  .   \end{array}$$
In the other hand,
$$\begin{array}{llll} 
\delta^n  H &  \equiv  \delta^{n-1}  H \delta + 2  \delta^n \mbox{ } 
\mbox{ } \mbox{ (By Proposition 4)}  \mbox{ } \\ & \equiv   \delta^{n-2} \; H \; 
\delta^2 +2 (2 \delta^n) \equiv  \cdots  \\ & \equiv
 n( 2 \delta^n) + H \delta^n (\rm{mod} \; \mathcal{U} \delta).
\end{array}$$ 
So we have: 
$$\begin{array}{llll}    
\delta^j \; D^j  & \equiv \delta D \; \delta^{j-1} \; D^{j-1} + 
\sum_{n=1}^{j-1} (H + 2n) \delta^{j-1} D^{j-1} \\  & \equiv
\delta D \; \delta^{j-1} \; D^{j-1} + (j-1)(H+j) \delta^{j-1} D^{j-1}
(\rm{mod} \; \mathcal{U} \delta),
\end{array}$$
Finally:
$$\begin{array}{llll} P_j(H)  & \equiv  H P_{j-1}(H) + (j-1)(H+j) P_{j-1}(H)  \\ 
 & \equiv  j(H+(j-1)) \; P_{j-1}(H)(\rm{mod} \; \mathcal{U} \delta).
\end{array}$$ We obtain the result: 
$P_n \equiv   n! \prod_{j=0}^{n-1} (H+j)(\rm{mod} \; \mathcal{U}
\delta) $ by induction.
\end{proof}
\begin{corollaire}
Let $f \in M_k$ a modular form of weight $k$ and $n \geq 0.$ We have
an exact sequence:
$$\delta^n \; D^n (f) =n!^2  \;  \binom{k+n-1} {n}  \; f.$$
\end{corollaire}

\begin{proof}
By using $(i)$ of   Proposition $2$ , $f \in \ker(\delta).$
So the last Proposition implies the result. $\Box$
\end{proof}
\begin{proposition}
Let $\Gamma \subset \rm{PSL}(2, \RM)$ be a discrete
and finite covolume subgroup. Let $k \geq 0$
and $p \geq 0$ be integers, if $p < \frac{k} {2}$ then:
$$\widetilde{M}_k^{(\leq p)} = D^p(M_{k - 2p}) \oplus
\widetilde{M}_k^{(\leq p - 1)}.$$ 
\end{proposition}
\begin{proof}
By Proposition $\ref{pro:del:pres}$ and its corollary we have, 
$\delta ^p(f) \in M_{k-2p}.$ By application of the corollary 
of Proposition $\ref{pro:cal:ope}$ to $\delta ^p(f)$ 
we get: 
$$p!^2\binom{k-p-1} {p}  
f - D^p( \delta^p(f)) \in \widetilde{M}_k ^{(\leq p-1)}(\Gamma).$$ 
In particular, if $k > 2p$ then $f$ is the
sum of the  $pth$ derivative of a  modular 
form and of a quasimodular form of depth $\leq p-1.$
\end{proof}
We finish this section by giving an additive structure
theorem and an $\rm{sl}_2(\CM)$-module structure theorem
for rings of quasimodular forms
over discrete and finite covolume subgroups 
of $\rm{PSL}(2, \RM).$ 

\begin{theoreme}
Let $\Gamma \subset \rm{PSL}(2, \RM)$ be a discrete
and finite covolume subgroup. We have an exact sequence:
$$0 \longrightarrow M_2(\Gamma) \longrightarrow 
\widetilde{M}_2(\Gamma) \stackrel{\delta}\longrightarrow \CM .$$
Then we have:
$$\widetilde{M}_* = \CM \oplus \bigoplus_{i=0}^{\infty}
(D^i (M_*)) \oplus  \left\{ 
\begin{array}{cccc}  
 0 \mbox{ } &\mbox{ } \mbox{if} \mbox{ }  \dim(\rm{Im}(\delta))  = 0. \\
\oplus_{i = 0}^{\infty}  \CM D^i ( \phi)   &\mbox{ } \mbox{if}  \mbox{  }  
\mbox{ there exists} \mbox{ } \phi \in \widetilde{M}_2(\Gamma) 
\mbox{  } \\     
& \mbox{  } \mbox{satisfying}  \mbox{ } \delta(\phi) = 1. 
\end{array} 
\right . $$
\end{theoreme}
\begin{proof}
By general properties of quasimodular forms, it is clear
that $M_2(\Gamma) \subset \widetilde{M}_2(\Gamma)$ and
$\delta(M_2(\Gamma)) = 0.$ 

We suppose that $\dim \rm{Im}(\delta) = 0.$ Then for any
$f \in \widetilde{M}_k^{\leq p},$ we have  $p < \frac{k} {2}.$ 
We deduce from Proposition $5$ that $f$ is the direct sum of
the $p$th derivative of a modular form and  of a quasimodular
form of depth $< p.$ So by induction on $p,$ we get that
$\widetilde{M}_k = \bigoplus_{i=0}^{\frac{k} {2}} D^i M_{k - 2i}.$
This implies the theorem in the case $\dim(\rm{Im}(\delta)) = 0.$

We suppose know that there exists $\phi$ such that $\delta(\phi)=1.$
By using the same proof for the additive structure of 
$\widetilde{M}_*(\Gamma_1)$ (as in $[4]$), we deduce the theorem
in the second case. 
\end{proof}
For $k > 0,$
let $\mathcal{U}_k$ be the $\rm{sl}_2(\CM)$-module defined
by a basis $(x_j^{(k)})_{j \in \NM}$ 
with $D x_j^{(k)} = x_{j+1}^{(k)},$ 
$H x_j^{(k)} = (k + 2j) x_j^{(k)}$ and 
$$\delta x_j^{(k)} = \left\{ 
\begin{array}{ccc}
j(k+j-1) x_{j-1}^{(k)} & \mbox{if} \mbox{ } j \geq 1. \\
0 & \mbox{if} \mbox{ } j = 0. \end{array} \right.$$
We define $\mathcal{U}_0 = \CM$ with the trivial action of
$\rm{sl}_2(\CM).$ Finally, if there exists $\phi$ such that
$\delta(\phi)=1,$ we define an extension $\widehat{U}_2$
of $\mathcal{U}_2$ (i.e. $\widehat{U}_2 \simeq  \CM \oplus
\mathcal{U}_2$). If  $(\hat{x}_j)_{j \in \NM}$ is a basis of
$\widehat{U}_2$ then the the action of $\rm{sl}_2(\CM)$
is defined as over  a basis of $\mathcal{U}_2$ 
except that $\delta \hat{x}_0 = 1.$  

For any $k > 0,$ we have an embedding:
$$\begin{array}{cccc}
M_k \otimes \mathcal{U}_k &  \longrightarrow & \widetilde{M}_* \\
f \otimes x_j^{(k)}  & \longrightarrow &  D^j f. \end{array}$$
In the case $k=0,$ we have a map $\CM \otimes \CM \longrightarrow
\CM.$ Finally 
$$\begin{array}{cccc}
\CM \phi \otimes \widehat{U}_2 & \longrightarrow & \widetilde{M}_* \\
\phi \otimes \hat{x}_i & \longrightarrow & D^i \phi \mbox{ } \mbox{if} \mbox{ }
i \geq 1 \\ \phi \otimes 1 & \longrightarrow & 1 \mbox{ } \mbox{if}
\mbox{ } i = 0
\end{array}$$
\begin{theoreme}
Let $\Gamma \subset \rm{PSL}(2, \RM)$ be a discrete and 
finite covolume subgroup. Then we have:
$$ \widetilde{M}_*(\Gamma) = \bigoplus_{k = 0}^{\infty}
M_k(\Gamma) \otimes \mathcal{U}_k \oplus \left\{ \begin{array}{ccc}
0  & \mbox{ } \mbox{if} \mbox{ }  \mbox{ } 
\delta(\widetilde{M}_2(\Gamma)) = 0. \\
\widehat{U}_2 & \mbox{ } \mbox{if} \mbox{ }  \mbox{ }
\mbox{there exists} \mbox{ }  \phi \in \widetilde{M}_2(\Gamma)  
\\  & \mbox{such that} \mbox{  } 
\delta(\phi) = 1. \end{array} \right.$$ 
\end{theoreme}
The proof of this theorem use
the definitions of the maps $M_k \otimes \mathcal{U}_k
\longrightarrow \widetilde{M}_*,$ and Theorem $2.$

\section{The cocompact / non-cocompact dichotomy}
We  prove a dichotomy Theorem which characterizes
cocompact groups in terms of their space of
quasimodular forms of weight $2.$ 
This fundamental dichotomy implies
the  condition $p < \frac{k} {2}$  in the case
of cocompact groups and we deduce from it the additive structure 
theorem in the cocompact case.

\begin{theoreme}
{\bf Let $\Gamma \subset \rm{PSL}(2,\RM)$ 
be a discrete and  finite covolume  subgroup. 
If  $\Gamma$ is not cocompact, there exists
a quasimodular form $\phi$ of weight $2$ over
$\Gamma;$ morover $\phi$ is not modular and
$\widetilde{M}_2 = M_2(\Gamma) \oplus \CM \phi.$
If $\Gamma$ is cocompact we have
$\widetilde{M}_2 (\Gamma) = M_2 (\Gamma).$}
\end{theoreme}
\begin{proof}
In the case of a  group $\Gamma$ which is commensurable with
$\Gamma_1 = \rm{PSL}(2, \ZM),$ we take $\phi$ equal
to the restriction of the Eisenstein series $E_2$
if $\Gamma$ is a subgroup of congruence of $\Gamma_1$
and in the other cases  of a group $\Gamma$ 
commensurable with $\Gamma_1,$  take a normalised trace of
$E_2$ and get a quasimodular form of weight $2$
over $\Gamma$ with $\delta(\phi)=1.$ In particular,
$\phi$ is not modular. If $\Gamma$ is not an arithmetic
group (such as a Hecke modular group) then  we can define
$\phi$ as the quasimodular form associated to an almost
holomorphic modular form $E_{2,\Gamma}(z,0)$ of weight $2$
defined as the limit in $s$ of a family
$E_{2, \Gamma}(z,s)$ of almost holomorphic
modular forms of weight $2.$

We suppose that $\Gamma$ is cocompact and  that there exists  a
quasimodular form $f$ of weight $2$ 
which is not modular. Let $F$ be the almost holomorphic modular form  
associated to $f.$
We have : $$ F(z) =f(z) + \frac{c} {z - \overline{z}} 
\mbox{  } \mbox{with}  \mbox{ } c \neq 0,$$ in fact, 
$f_0 = f \in M_2 $ so $f_1 \in M_0 = \CM.$
Let $\omega (z)= F(z)\; dz,$ the modularity of  $F$ 
implies the $\Gamma$ invariance of the  $\omega$ form. 
So this $1-$form is defined on the quotient $X = 
\mathcal{H} / \Gamma.$ On the other hand, we have  $$d 
\omega = - \frac{\partial F} {\partial \overline{z}} \;  
dz \wedge d \overline{z} = - \frac{c} {(z - \overline{z})^2} 
\; dz \wedge d\overline{z}.$$ 
This means that $d \omega$ is a multiple of the volume form. So
there exists $\alpha \neq 0$ such that : 
$$0 \neq \alpha \; \rm{Vol}(X) = \int_X d \omega.$$ 
On the other hand $ \int_X d \omega  = 0,$ this equality
is a consequence of Stokes's formula and the fact that $X$
is a variety without boundary. We obtain a contradiction. 
\end{proof}
\begin{remarque}
Theorem $4$ implies that $\delta(\widetilde{M}_2) =
\delta(M_2) = 0$ in the cocompact case. By Theorem $2$
we deduce the corresponding additive structure for
rings of quasi-modular forms. In the non-cocompact
case $\delta(\widetilde{M}_2) = \CM$ and by Theorem
$2$, we deduce from this the corresponding additive structure.
\end{remarque}
\section{Rings of quasimodular forms}
%%%%%%%%%%%%%%%%%%%%%%%%%%%%%%%%%%%%%%%%%%%%%%%%%%%%%%%%%%%%%%%%%%%%%%%%
%%%%%%%%%%%%%%%%%%%%%%%%%%%%%%%%%%%%%%%%%%%%%%%%%%%%%%%%%%%%%%%%%%%%%%%%
For a discrete and cocompact subgroup $\Gamma$ we denote
by  $I$ (respectively $\widetilde{I}$) the ideal of modular
forms (respectively quasimodular forms) over $\Gamma$ of 
positive weights.
Finally, $\widetilde{I}^2 _k = \sum\limits_{0 < j < k} 
\widetilde{M}_j  \widetilde{M}_{k-j}$ is the  
$\CM$-vector space of decomposable quasimodular forms 
of weight $k.$ 

\begin{theoreme}
{\bf Let $\Gamma \subset \mathrm{PSL(2,\mathbb{R})}$ be a 
discrete and  cocompact subgroup. Let $\epsilon = \dim_{\CM}  
I/(I \cap  \widetilde{I}^2)$ and let $\{ A_1 , \cdots , A_{\epsilon} \}$ 
be homogeneous elements of $I$ of weights  $w_1, \cdots ,w_{\epsilon} 
\in 2 \ZM,$ which are  linearly independent modulo
$(\widetilde{I})^2$. Then for any $k \geq 0,$  
$$ (\widetilde{I}/ \widetilde{I}^2)_k = 
\bigoplus_{i=1,w_i \leq k}^{\epsilon} \mathbb{C} D^{( \frac{k - 
w_i} {2})} (A_i) \quad\quad.$$} 
\end{theoreme}
\begin{proof}
We denote by $P_s,$ ($s = 2,4,\cdots$) the vector space generated
by all the  $A_i$  of weight $w_i =s$ and we write
$\delta_i = \dim P_i$ so $ \sum_i \delta_i = \epsilon .$ 
We have the commutative diagram: 
$$\begin{array}{cccc}   P_s & \hookrightarrow  &  M_s \\  & \searrow & 
\downarrow \quad\llap{D}^n \\ &  &  \widetilde{M}_{s +2n}\end{array} $$
The maps which appear in the last diagram are injective. In fact
$P_s \subset M_s$. On the other hand  $D^n(f) = 0$ implies that $f$ is 
a polynomial. This polynomial is zero because it  defines a modular
form $f$ of positive weight. We deduce: 
$$\dim D^n (P_s) = \delta_s\quad.$$  On the other hand 
$D^{\frac{k -2} {2}}(P_2) \subset   D^{\frac{k - 2} {2}}(M_{2}), 
\cdots, D^{\frac{k - w_{\epsilon}} {2}}(P_{w_{\epsilon}}) 
\subset   D^{\frac{k - w_{\epsilon}} {2}}(M_{w_{\epsilon}}).$
Since Theorem $2$  implies that the sum of spaces 
$D^{\frac{k-s} {2}} M_s,$  for $s=2,4,\cdots ,w_{\epsilon}$, is a direct 
sum, we get that the sum of their subspaces $D^{\frac{k - s} {4}}
P_s$ is also  direct.
On the other hand, for any $n\geq 0$ and  $s$ : 
$2 \leq s \leq w_{\epsilon}$ 
we have, $$ D^n(P_s) \cap (\widetilde{I})^2 = 0 \quad\quad.$$  
Indeeed, by corollary of Proposition $\ref{pro:cal:ope},$ 
$$\mbox{for all} \mbox{ } f \in P_s \mbox{ } \mbox{,} \mbox{ } \delta ^n D^n (f) 
= c_n \; f \mbox{ } \mbox{with} \mbox{ } c_n \neq 0 \quad.$$ Since 
$\delta$ is a derivation, we have: $$ \mbox{for all} \mbox{ } f , g \in 
\widetilde{I}, \delta(g h)=\delta(g) h + g \delta(h)\quad\quad.$$ 
This implies $ \delta( \widetilde{I}^2) \subset \widetilde{I}^2,$ 
indeed   $M_0 \cap Im(\delta) = 0 $ because $\widetilde{M}_2 = M_2.$ 
It remains to prove that : for any $ f_2 \in P_2 , \cdots , 
f_{w_{\epsilon}} \in P_{w_{\epsilon}},$ 
if $f_2^{( \frac{k - 2} {2}) } + \cdots + f_{w_{\epsilon}} ^{(
\frac{k-w_{\epsilon}} {2})}
\in \widetilde{I}^2$ then $ f_2 = \cdots = f_{w_{\epsilon}} = 0. $ 
We write $\alpha_2 = \frac{k - 2} {2}, \cdots ,\alpha_{w_{\epsilon}}
= \frac{k - w_{\epsilon}} {2}$ and 
we suppose that: $$f_2 ^{(\alpha_2)} + \cdots + f_{w_{\epsilon}} ^
{(\alpha_{\epsilon})} \in \widetilde{I}^2,$$ 
we can suppose also that:  $$\alpha_2 \geq \alpha_4 \geq \cdots \geq 
\alpha_{w_{\epsilon}}.$$ 
We apply the operator $\delta^{\alpha_2},$ then all 
$ f_i ^{ (\alpha_i)}$ with $i > 2$ vanish. We deduce from this that   
$f_2 \in \delta ^{(\alpha_2)} (\widetilde{I}^2) \subset \widetilde{I}^2$ 
and so $f_2=0.$ We restart with the operator 
$\delta^{\alpha_4},$ we prove that $f_4=0$ and by induction,  
we deduce that $f_{w_{\epsilon}} = 0.$ 
\end{proof}
\begin{corollaire}
{\bf Let $\Gamma \subset \mathrm{PSL(2,\mathbb{R})}$ be a 
discrete and  cocompact subgroup. Let $\epsilon = 
\dim_{\mathbb{C}}  I / ( I \cap \widetilde{I}^2)$ and  let $\{ A_1 , \cdots ,
A_{\epsilon} \}$ be  homogeneous elements of $I,$ linearly independent modulo 
$\widetilde{I}^2 $ of  weights $w_1, \cdots ,w_{\epsilon}$ respectively, 
then: $$\dim_{\mathbb{C}} ( \widetilde{I} / \widetilde{I} ^2)_k = 
\epsilon, \mbox{ } \mbox{for all} \mbox{ } k \geq \max_i \{w_i\}$$ 
In particular, $\widetilde{M}_*$  is a non-finitely generated $\CM$-algebra.}
\end{corollaire}
\begin{remarque}
{\rm This result is false in the  non-cocompact case. 
For example for $\rm{PSL}(2,\mathbb{Z}),$ we have
$ \widetilde{M}_* \simeq \mathbb{C}[E_2,E_4,E_6]$ where $E_2$,
$E_4$ and  $E_6$ are the Eisenstein series of weights 
$2$,$4$ and $6,$  respectively.}
\end{remarque}
\begin{proof}
The fact that $\widetilde{M}_*$ is not finitely generated 
is equivalent to saying that
$\dim((\widetilde{I} / \widetilde{I} ^2 )_k) =0,$ for $k$ 
large enough. The corollary is a consequence of Theorem $~3.$ 
\end{proof}
\begin{remarque}
We finish this section by observing that the ideas 
used to prove our two theorems can be used to prove more 
general results about the additive and  multiplicative 
structure of the differential closure $\rm{CL}(\mathcal{M})
_*$ of any ring $\mathcal{M}_*$ different from $\mathbb{C}$ 
and generated by a finite set of 
holomorphic modular forms  or meromorphic 
modular forms of  positive weight. We will define the 
differential closure, and we give the corresponding results
(Theorems $6$ and $7$). The main point for Theorem $7$
is the fact that $\rm{CL}(\mathcal{M})_2 = \mathcal{M}_2.$
\end{remarque}
\begin{definition}
Let $\Gamma \subset \rm{PSL}(2, \mathbb{R})$ be a discrete 
subgroup of finite covolume, and let $\mathcal{M}_*$ be a graded 
subring of the ring of meromorphic modular forms over $\Gamma$.
The differential closure $\rm{CL}(\mathcal{M})_*$ of $\mathcal{M}_*$ 
is the  smallest ring which contains $\mathcal{M}_*$ and is closed 
under the derivation $D,$ with  graduation 
$D^j \mathcal{M}_k \subset \rm{CL}(\mathcal{M})_{k + 2j}.$
\end{definition}

\begin{notation}
We denote by $J_{\mathcal{M}}$ the ideal 
of elements in $\rm{CL}(\mathcal{M})_*$
of positive weight and by $I_{\mathcal{M}}$ the ideal of 
$\mathcal{M}_*$ of 
elements of positive weight.  The ideal $J_{\mathcal{M}} ^2$ is 
the ideal of  decomposable forms in $\rm{CL}(\mathcal{M})_*.$ 
\end{notation}

\begin{theoreme}
{\bf Let $\Gamma \subset \rm{PSL}(2,\mathbb{R})$ be a 
discrete subgroup of finite covolume, and let $\mathcal{M}_*$ be   
a ring generated by a finite set of  holomorphic modular forms 
or meromorphic modular forms  of  positive weight over $\Gamma.$ 
Then for any $k \geq 0$:
$$ \rm{CL}(\mathcal{M})_k =\bigoplus_{0 \leq j \leq \frac{k} {2} } 
D^j \mathcal{M}_{k-2j}.$$}
\end{theoreme}
\begin{theoreme}
{\bf Let $\mathcal{M}_*$ be a ring like the one of the last theorem. 
Let $\epsilon = \dim_{\CM} I_{\mathcal{M}}/(I_{\mathcal{M}} 
\cap J_{\mathcal{M}}^2)$ 
and let $\{f_1, \cdots, f_{\epsilon} \}$ be homogeneous elements
of  $I_{\mathcal{M}},$ linearly independent
modulo $J_{\mathcal{M}}^2 $ of weights 
$l_1,\cdots,l_{\epsilon},$ respectively. Then for any $k$ even:  
$$ (J_{\mathcal{M}} / J_{\mathcal{M}} ^2)_k = 
\bigoplus_{\begin{array}{cc}  i=1,\cdots,\epsilon \\
l_i \leq k \end{array}} \CM D^{(\frac{k-l_i} {2})}(f_i).$$
In particular, 
$$\dim_{\CM}(J_{\mathcal{M}} / J_{\mathcal{M}}^2)_k = \epsilon, 
\quad\quad \mbox{for all} \mbox{ } k \geq \max \{l_1,\cdots,l_{\epsilon} \}\,,$$
and the ring $\rm{CL}(\mathcal{M})_*$ is not finitely generated as a 
$\CM$ algebra.}
\end{theoreme}
%%%%%%%%%%%%%%%%%%%%%%%%%%%%%%%%%%%%%%%%%%%%%%%%%%%%%%%%%%%%%%%%%%%%%%%%%
%%%%%%%%%%%%%%%%%%%%%%%%%%%%%%%%%%%%%%%%%%%%%%%%%%%%%%%%%%%%%%%%%%%%%%%%%
\section{Algebraic characterization of cocompact modular groups}
We recall that a Poisson algebra is    
a commutative and  associative algebra $A$
with a Lie structure, i.e. a bilinear operation 
$[\,\cdot \, , \, \cdot \,] : A \times A \longrightarrow A$ 
satisfaying the Jacobi identity, such that for any 
$x \in A,$ the map $[x,\,\cdot\,]$ is a derivation. 
If furthermore $A= \bigoplus_{n \geq 0}
A_n$ is graded with $A_m A_n \subset A_{m+n},\; [A_m,A_n] \subset 
A_{m+n+1},$ then $A$ is called a graded Poisson algebra.

\noindent{\it Examples. }
1)\;Let $A$ be a graded algebra (commutative and associative) 
and let $d: A \longrightarrow A$ be  a  derivation of degree $1,$ i.e. 
$d(A_n) \subset A_{n+1}$ and $d(xy) = x d(y) + y d(x),$ 
for evry $x,y \in A.$ Then the braket
defined by $[x,y] = H(x) d(y) - H(y) d(x),$ where $H$ is 
the operator of  multiplication by the  
weight $n$ in $A_n$ satisfies the Jacobi identity 
(a simple verification) and has the property that:
$x \longrightarrow [x,y]$ is a derivation 
for evry fixed $y \in A$  (because $H$ and $d$ are derivations). 
We call a Poisson algebra  {\it trivialisable}, 
if it can be obtained in this way.

2)\;Let $\Gamma \subset \rm{PSL}(2, \RM)$ be a discrete
and finite covolume subgroup and let $A = M_{ev}= 
\bigoplus_{n \geq 0} M_{2n}$ (i.e., $A_n = M_{2n}$) be 
the graded algebra of modular forms. 
This algebra has a Poisson structure with the usual
multiplication and where the braket
$[\,\cdot\,,\,\cdot\, ]=[\,\cdot\,,\,\cdot\,]_1$, 
is the first Rankin-Cohen braket.

\begin{theoreme}
{\bf Let $\Gamma \subset \rm{PSL}(2, \RM)$ be a discrete
and finite covolume subgroup. Then 
the Poisson algebra $(M_{ev}(\Gamma) , [\,\cdot \, , \, \cdot \,]_1)$ 
is trivialisable if and only if $\Gamma$ is not cocompact.}
\end{theoreme} 
We use the  next two lemmas, where the second
is a  corollary of the first.
\begin{lemme}
Let $M_*^{\mathrm{mer}}$ be the ring of meromorphic modular 
forms over a discrete and finite covolume group,
then any derivation  $ \partial : M^{\mathrm{mer}}_* \longrightarrow 
M^{\mathrm{mer}}_{*+2} $ trivialising
the first Rankin-Cohen braket 
has the  form  $\partial = D - \phi E$ where $E$ is the
Euler operator (multiplication by the weight) and 
$\phi \in \widetilde{M}_2 ^{\mathrm{mer}}$ with $\delta \phi =1.$
\end{lemme}

\begin{proof}
For any discrete finite covolume group $\Gamma$ 
there exists a  meromorphic quasimodular  form 
$\psi$ of weight $2$ such that $\delta \psi =1$: we can divide
the logarithmic derivative of any non zero modular form
by its weight. 
Then, we consider $\partial_{\psi} =D -\psi E;$ 
this operator trivialises the first braket.
We suppose also that $\partial$ trivialises the first Rankin-Cohen braket. 
Then for any $f \in M^{mer}_k$ and $g \in M^{mer}_l$, 
we have the  relation:
$$\frac{\partial f - \partial_{\psi} f} { E(f)} = \frac{\partial g - 
\partial_{\psi} g} {E(g)}\quad\quad.$$
In particular, there exists a  holomorphic  modular form $\alpha$
of weight $2$ such that
$\partial - \partial_{\psi} = \alpha E$. This implies  
$\partial = D - (\alpha + \psi)E.$
Since $\psi + \alpha$ is a meromorphic quasimodular form 
of weight $2$ 
and  $\delta (\psi + \alpha ) =\delta \psi =1$ (because $\alpha$ 
is modular), we can take $\phi = \psi + \alpha$.
\end{proof}
\begin{lemme}
Let $M_*$ be the ring of modular forms over a discrete and finite 
covolume group $\Gamma$. 
Then any derivation 
$\partial : M_* \longrightarrow M_{*+2}$ trivialising 
the first Rankin-Cohen braket has the
form $\partial_{\phi} = D - \phi E$ with $\phi \in \widetilde{M}_2$ 
and  $\delta \phi = 1$.
\end{lemme}
\begin{proof}
By the last  lemma, $\partial$ has the form $\partial_{\phi}$ 
with $\phi \in M_2 ^{\rm{mer}}(\Gamma).$ 
The fact that $\partial_{\phi}(f)$ must be holomorphic for  
any holomorphic  modular form $f$ implies that  $\phi$ 
must  also be  holomorphic, 
since different modular forms over $\Gamma$ 
cannot have the  same set of zeros.
\end{proof}
The proof of Theorem $6$ is a   consequence of 
corollary of  Proposition $3$ and of the last lemma, 
because  holomorphic quasimodular forms
of weight $2$ over a cocompact modular 
group  are  modular (see Theorem$4$).
\begin{remarque}
This theorem is false in the non-cocompact case. 
for example, if 
$\Gamma = \mathrm{PSL}(2,\ZM)$ is the classical modular group; 
it then exists a derivation $\partial = D - \frac{E_2} {12} E$, 
where $E_2$ is the Eisenstein series of weight $2,$ 
which  trivialise the first Rankin-Cohen braket. 
\end{remarque}
%%%%%%%%%%%%%%%%%%%%%%%%%%%%%%%%%%%%%%%%%%%%%%%%%%%%%%%%%%%%%%%%%%%%%%%%%%%%%
%%%%%%%%%%%%%%%%%%%%%%%%%%%%%%%%%%%%%%%%%%%%%%%%%%%%%%%%%%%%%%%%%%%%%%%%%%%%%
\section{Emmbedding of quasimodular forms in finitely generated rings}
\begin{theoreme}
{\bf Let $\Gamma \subset \rm{PSL}(2, \RM)$ be a  discrete and  
cocompact subgroup.  Then there exists  a quasimodular 
form $\phi$ of weight $2$ over $\Gamma$ satisfying $\delta(\phi) = 1,$
with simple poles in the orbit of $i,$ and  without other
poles. For any such form $\phi$ we have 
${\rm Res}_{z=i} (\phi(z) dz)= \mathcal{K}$ 
for any $\alpha$ in the $\Gamma$-orbit of $i$ with
$\mathcal{K} = \frac{\rm{Vol}(\mathcal{H }/ \Gamma )} {4 \pi}$.}
\end{theoreme}
\begin{remarque}
{\rm The form $\phi$ is unique up to the addition of an holomorphic 
modular form of weight $2$  (the dimension of the space of such
forms is equal to the genus $g$ of the Riemann surface 
$\mathcal{H} / \Gamma$). }

After conjugating $\Gamma$ in $\rm{PSL}(2, \RM)$, we can replace 
$''i''$ in the  theorem by any other point
$z_0 \in \mathcal{H}.$
\end{remarque}
\begin{proof}
First, we suppose that $\Gamma$ acts on $\mathcal{H}$ without
fixed points (this means that the action is free).
Let $f$ be a non zero modular form of weight $k > 0,$  
we know that $\frac{f'} {f}$ is a meromorphic quasimodular form
of weight $2$ with $\delta(\frac{f'} {f}) = k \neq 0.$ Moreover
the  poles of $\frac{f'} {f} $ are simple 
and  $\Gamma$ invariant.
We denote by  $\{ P_1, \cdots , P_n \}$ the poles 
of  $\frac{f'}  {f}$  in $ \mathcal{H} / \Gamma$ different from $i.$
We want to construct a meromorphic modular form $h$ of weight $2$ 
such that the sum $ \frac{f'} {f}  + h$ has no poles outside 
the orbit of $i$.  Let $X = \mathcal{H} / \Gamma$ the 
Riemann compact surface (of genus $g$).
Then the hypothesis on $\Gamma$ implies that $X$ is smooth and that $g > 1.$
We denote by $\Omega^{1}_X$ the sheaf of holomorphic differential
$1$-forms over $X$. For any set of distinct points $\{q_1, 
\cdots, q_m\} \subset X$ 
(with $ m \geq 1$), we denote by $\Omega^{1} _X ( q_1 + \cdots + q_m)$ 
the sheaf of holomorphic differential $1-$ forms over $X$ with simple
poles at  $q_1,\cdots , q_m.$ We will prove: $$ H^0 (X ,
\Omega^1 _X (q_1 +\cdots + q_m))  \simeq   \mathbb{C}^{g +m-1}.$$
Let $K$ be the canonical divisor of $X.$ By the
Riemann-Roch Theorem we have: $$l(K + q_1 + \cdots +
q_m)=l( -(q_1 + \cdots + q_m)) + \deg( K + q_1 + \cdots +
q_m) -g +1.  $$ From $\deg(K) = 2 g -2 $ \; and \;  $l(-(q_1 + \cdots +
q_m))=0,$ we deduce that~: $$ l( K + q_1 + \cdots + q_m) = g + m
-1\quad. $$ If we apply the   Riemann Roch Theorem to the
cases $m=1$ and  $m=n+1,$ we obtain the exact sequence: 
$$ 0 \longrightarrow H^0(X,\Omega^1_X( i))
\longrightarrow H^0(X, \Omega^1_X( i + P_1 + \cdots + P_n))
\stackrel{\mathrm{Res}}\longrightarrow  \mathbb{C}^n 
\longrightarrow 0,$$
where $\rm{Res}$ maps 
a differential $1$-form $\omega$ to $(\rm{Res}_{P_1}(\omega), 
\cdots,\rm{Res}_{P_n}(\omega)).$
So we can choose $h$ of
weight $2$ such that $\phi = \frac{1} {k} \frac{f'} {f} + h $ has 
a simple pole at $i$ and no poles outside the orbit of $i.$ 
We also have $\delta \phi = 1.$

To compute the constant  $\mathcal{K}$ we apply Stokes's formula
to the  meromorphic  differential $1-$form   $\omega(t) =
\phi^*(t+i) \;  dt$  over $X$ where $\phi^*(t+i) = \phi(t+i) + 
\frac{1} {t - \overline{t} }$ is the almost holomorphic modular 
form associated to $\phi$.  Let $U_{\epsilon}$ be  a disk with center
$i$ and radius $\epsilon$  included in $X$. 
By Stokes's formula, we have: $$ \int_{X - U_{\epsilon}} d \omega(t) =
\int_{\partial (X - U_{\epsilon})} \omega(t)\quad$$ since, 
$d \omega(t) = d \phi^*(t)  \wedge dt =  - \frac{\partial \phi^* (t+i )} 
{\partial \overline{t} } d t \wedge d \overline{t}.$  On the other hand
$\phi$ is holomorphic over $\mathcal{H}$ so: 
$$- \frac{ \partial \phi^* } {\partial \overline{t}} =  
\frac{\partial} {\partial \overline{t}} ( \frac{1} {t - \overline{t}}) 
\,dt \wedge d \overline{t}. $$ because $\phi$  satisfies $ \frac{\partial}
{\partial \overline{t} } \phi = 0.$ We obtain $d \omega (t ) =
\ \frac{ dt \wedge d \overline{t} } { ( t - \overline{t} )^2 },$ 
or $\frac{1} {2i}$ times
the volume form, so $\int_{X - U_{\epsilon} } d \omega(t) = \frac{1} {2i} 
\rm{Vol}( X - U_{\epsilon}).$ On the other hand 
$\int_{\partial (X - U_{\epsilon} )}
\omega(t) = - \int_{\partial U_{\epsilon} } \omega (t) $ because $X$
is a compact manifold without boundary. So $\int_{\partial (X -
U_{\epsilon} )} \omega(t) = -\int_{\partial(U_{\epsilon})}
\phi(t+i) + O(1),$ indeed $\phi^*(t+i) - \phi(t+i)$ is a continuous map 
over $\partial U_{\epsilon}$. In the other hand $\phi(t+i) \sim
\frac{\mathcal{K}} {t}$ so $ - \int_{\partial U_{\epsilon}}
\omega(t) = - (2 \pi i) \mathcal{K} + O(1);$ by letting
$\epsilon$ to $0$ we obtain $\mathcal{K} =  \frac{ \rm{Vol}(X)} 
{4 \pi}.$  This completes the proof in the case  of
groups acting on  $\mathcal{H}$ without fixed points.

We suppose now  that $\Gamma$ acts on  $\mathcal{H}$ and that the action
is not necessarily free.  Selberg lemma implies that there exists a 
subgroup $\Gamma ' \subset \Gamma $ of  finite index  without
torsion. The first part of the proof implies that there exists 
a quasimodular form $\alpha$ over $\Gamma'$ of weight $2$ with at most simple
poles in the orbit of the point $i$.
We define: $$ \beta (z) =\sum_{ \gamma \in
\Gamma / \Gamma'} [(\alpha \mid \gamma)(z) - \frac{c} {c z + d}],$$
with $\gamma = \left( \begin{array}{cc} a & b \\ c & d \end{array}
\right)$ and $( \alpha \mid \gamma )(z) =(c z + d)^{-2} \alpha (\frac{
a  z + b} {c  z + d} ).$ We will prove that $\beta$ is a 
quasimodular form over $\Gamma$ of weight $2$. Let $\alpha^*$ 
be the almost holomorphic modular form associated to $\alpha.$ 
It is easy to check that
$$ \beta^*(z) = \sum_{ \gamma \in  \Gamma / \Gamma' } (\alpha^*
\mid \gamma)(z), $$ is an almost holomorphic modular form over $\Gamma$ of
weight $2$ (since $\alpha^*$ is modular, $\beta^*$ corresponds 
to the trace of $\alpha^*$ over the group $\Gamma$). On the other hand 
we have 
$(\alpha^* \mid \gamma) (z) = [ ( \alpha \mid \gamma) (z) - \frac{c} {
c z + d}  ]  + \frac{1} { z - \overline{z} }.$ So: $$ \beta^*(z) =
\sum_{ \gamma \in \Gamma / \Gamma '} [ ( \alpha \mid \gamma )(z) -
\frac{c} { c z + d} ] + \sum_{ \gamma \in \Gamma / \Gamma'} \frac{1}
{ z - \overline{z}},$$ in other words $\beta^*(z) = \beta(z) + \frac{
[\Gamma : \Gamma']} { z - \overline{z} }.$ This proves that $\beta$
is a  quasimodular form over $\Gamma$ of weight $2$ and $\delta (\beta ) = 
[\Gamma:\Gamma'] $. It is clear that $\beta$ has at most 
simple poles on the  orbite of  $i$. Hence,  $\frac{\beta}
{[\Gamma:\Gamma']}$ is an appropriate form over $\Gamma$.
\end{proof}
\begin{notation}
We denote by $M_*(\Gamma;\{i\})$ the ring of modular forms 
without poles outside
the orbit of $i$ and we denote by  $M_* ^{ (\geq
\alpha) } (\Gamma ; i)$ the subset of modular forms  over
$\Gamma $ with vanishing order at least equal to $\alpha$ 
at $i$. Finally we denote by
$\widetilde{M}_2(\Gamma;\{ i \}),$ the space of quasimodular forms 
of weight $2$ over $\Gamma$ with all poles in the orbit of  $i.$
\end{notation}
\begin{lemme}
Let $\Gamma \subset \rm{PSL}(2, \RM)$ be a discrete cocompact subgroup,
and  $\phi$ a quasimodular form over $\Gamma$ with at most
simple poles in the orbit of $i$, and
$\delta(\phi)=1.$ Then, we have ${\rm Res}_i(\phi(z) dz)=
\frac{{\rm Vol}(\mathcal{H} / \Gamma)} {4 \pi}$ and 
$\omega = \phi' -\phi^2$
is a modular form of weight $4$ with at most double poles
in the orbit of $i.$
\end{lemme}
\begin{proof}
We know that for any 
$\left(\begin{array}{cc} a & b \\ c & d\end{array} \right)  \in \Gamma$ 
we have: $$ \phi( \frac{ a z +b} { c z +d} )= (c z +d)^2 \phi( z ) +  
c \; ( c z +d)\quad\quad.$$ By derivation, we  get: 
$$ \phi'(\frac{ a z + b } { c z + d} ) = (c z + d)^4 \phi'(z )  + 2  
c (c z + d)^3 \phi (z) + c^2 (c z + d)^2.$$ On the other hand: $$
\phi^2(\frac{ a z + b} { c z + d} ) = ( c z + d)^4 \phi^2(z) + 2 
c (c z +d)^3 \phi(z) +  c^2 (c z + d)^2, $$ this implies: 
$$(\phi' -  \phi^2) (\frac{ a z + b} { c z + d}) = 
(c z + d)^4 (  \phi' -  \phi^2) (z).$$ So $\omega$ is
a modular form of weight $4.$ Since $\phi'(i+x)  \sim
-\mathcal{K}  x^{-2}$ and  $\phi^2(i + x) \sim \mathcal{K}^2 x^{-2}$ 
(for $x \rightarrow 0$), we deduce that: 
$$\omega(x+i)  \sim -\mathcal{K} (\mathcal{K} + 1) x^2.$$
\end{proof}
\begin{proposition} 
Let  $\Gamma \subset \rm{PSL}(2, \RM)$ be a discrete cocompact 
subgroup, and $\phi$ a quasimodular form of weight $2$ over
$\Gamma$ with $\delta(\phi)= 1$, which is holomorphic  outside 
the orbit of $i.$
Then there exists an operator $$ D_{\phi}: M_k( \Gamma; \{ i \}) 
\longrightarrow M_{k +2}(\Gamma; \{ i \}),$$
defined by  $D_{\phi}(f) = f' - k \phi f.$ 
If $\phi$  has a simple pole at $i,$ then
$\rm{ord}_i(D_{\phi}(f)) \geq
\rm{ord}_i(f) -1,$ with inequality if and only if  
$\rm{ord}_i(f) = k \, \mathcal{K},$ where $k$ is the weight of 
$f.$
\end{proposition}
\begin{remarque}
{\rm The case $\rm{ord}_i(D_{\phi} (f)) = \infty$, can
only happen if $\rm{ord}_i(f) =  k(f) \mathcal{K}$ where $k(f)$
is the weight of $f$.}
\end{remarque}
\begin{proof} 
Let $f$ be a meromorphic modular form  over $\Gamma$ of weight $k$,
then  for any $\left( \begin{array}{cc} a & b \\ c & d
\end{array} \right)  \in \Gamma,$ we have  $ f(\frac{ a z +b } { c z + d} )
=(c z +d )^k f(z).$ By derivation, we obtain:$$ D f( \frac{ a z +b }
{ c z + d} ) = ( c z + d )^{ k + 2 } Df (z) + k c ( c  z + d) ^{ k + 1
} f(z). $$ On the other hand:  $$  ( \phi .f )( \frac{ a z + b} { c z + d} )
= ( c z + d) ^{k+2} ( \phi . f) ( z) + c ( c z +d)^{k+1} \; f(z).$$ This
implies  $D_{\phi}(f) (\frac{a z + b} {c z + d}) = (c z + d)^{k+2}
D_{\phi} (f)(z).$ On other words $D_{\phi}(f)$ is a modular form
of weight $k+2.$ In the other hand, if $f(x) \sim  x^{\alpha}$
(with $\alpha \neq \mathcal{K} k $) then  $D_{\phi}(f)(x) \sim (\alpha - 
k \; \mathcal{K} ) \; x^{\alpha -1}.$
\end{proof}
We recall that $I$ is the ideal of modular forms 
of positive weight over the group $\Gamma.$ 
\begin{theoreme}
{\bf Let $\Gamma \subset \rm{PSL}(2, \RM)$ be a discrete  cocompact
subgroup,
let $\phi \in \widetilde{M}_2 (\Gamma ;\{i \})$ with $\delta ( \phi ) =1$ 
and $\omega =\phi' -   \phi^2.$  Let $(f_1, \cdots ,f_d)$ be a basis
of $I/I^2$. Then there exists
$N \in \mathbb{N}^*$ such that the ring $R$ generated by 
the set: $$\{ D_{\phi} ^j (f_i) \; (1 \leq j \leq N , 1 \leq i \leq d) \; ;  
D_{\phi} ^l (\omega)\; (1 \leq l \leq N) \},$$ is closed 
under the operator $D_{\phi}.$}
\end{theoreme} 
To prove  the  theorem, we use a lemma which 
describe finitely generated semigroups of  
$\mathbb{R}^2$:
\begin{lemme}
Let  $G$ be a  finitely generated semigroup of $\RM^2$.  
We suppose that the group generated by $G$ is  
a lattice $\Lambda \subset \mathbb{R}^2$ of rank $2$. 
Let $S$ be the  sector $<G \,\cdot\, \mathbb{R}_{+}>.$  
We suppose that $S$ is  convex, with angle at most equal 
to $\pi.$ Then there exists  $A \in S$ such that 
$(A+S) \cap \Lambda \subset G.$
\end{lemme}
\begin{proof}
Let $\{P_1,\cdots , P_m \}$ be a  system  of generators of $G$. 
We suppose that the lines
$(O P_{m-1})$ and $(O P_m)$ bound the sector $S$.  
We consider a   coordinates system in 
$\RM^2$ in which $P_{m-1} =(1,0)$ and $P_m=(0,1)$.  
Then $\Lambda \otimes \QM = \QM^2 $ and 
the coordinates of  each $P_i$ are positive and rational
because $\{ P_{m-1},P_m \}$ is a  basis of 
$\Lambda \otimes_{\ZM} \QM $ over $\QM$.  
Let $P =(x,y) \in S \cap \Lambda$ be any point, then there
exists $(\alpha_1,\cdots ,\alpha_m) \in \mathbb{Z}^m$ such that:
$$ P = \alpha_1 P_1 + \cdots + \alpha_m P_m. $$ 
For any $i= 1, \cdots , m -2,$ there exists  
$ \overline{\alpha_i} $,\; $ 0 \leq \overline{\alpha_i} < a_i $ such that: 
$ \alpha_i \equiv \overline{\alpha_i}\,(\!\mod a_i).$ So, we can write 
$P$ as ~: $$ P= \overline{ \alpha_1} P_1 + \cdots +  
\overline{\alpha_{m_2}} P_{m-2} + 
\beta P_{m-1} + \gamma P_m,$$ with $\beta,\,\gamma \in \ZM.$
If the abscissa of $P$ satisfies: $$x(P) \geq X_0 := \max_{\{\; 0 
\leq \overline{\alpha_1} < a_1 ,\; \cdots\; , 0 \;\leq \overline{\alpha_{m-2}} 
< a_{m-2} \}} \;  x( \overline{\alpha}_1 P_1 + \cdots + 
\overline{\alpha}_{m-2} P_{m-2} ) $$ 
then $\beta \geq 0$. If the
ordinate of  $P$ satisfies: $$ y(P) \geq Y_0 := \max_{ \{ \; 0 \leq 
\overline{\alpha_1} \leq  a_1,\; \cdots \; , 0 \; \leq \overline{\alpha_{m-2}} 
\leq a_{m-2} \; \} } 
\; y(\overline{\alpha}_1 P_1 + \cdots + \overline{\alpha}_{m-2} P_{m-2}) $$ 
then $\gamma \geq 0$. We  take  $A=(X_0,Y_0).$ 
\end{proof}
Back to the proof of Theorem $7.$
\begin{proof}
We consider the map:$\!$
$$\begin{array}{ccc} I: M_*(\Gamma ;
\{ i \})  & \longrightarrow & \mathbb{N}^2 \\ f & \longrightarrow
&(\frac{k(f)} {2} , \rm{ord}_i(f) + \frac{k(f)} {2}),
\end{array}$$
\begin{figure}
\includegraphics{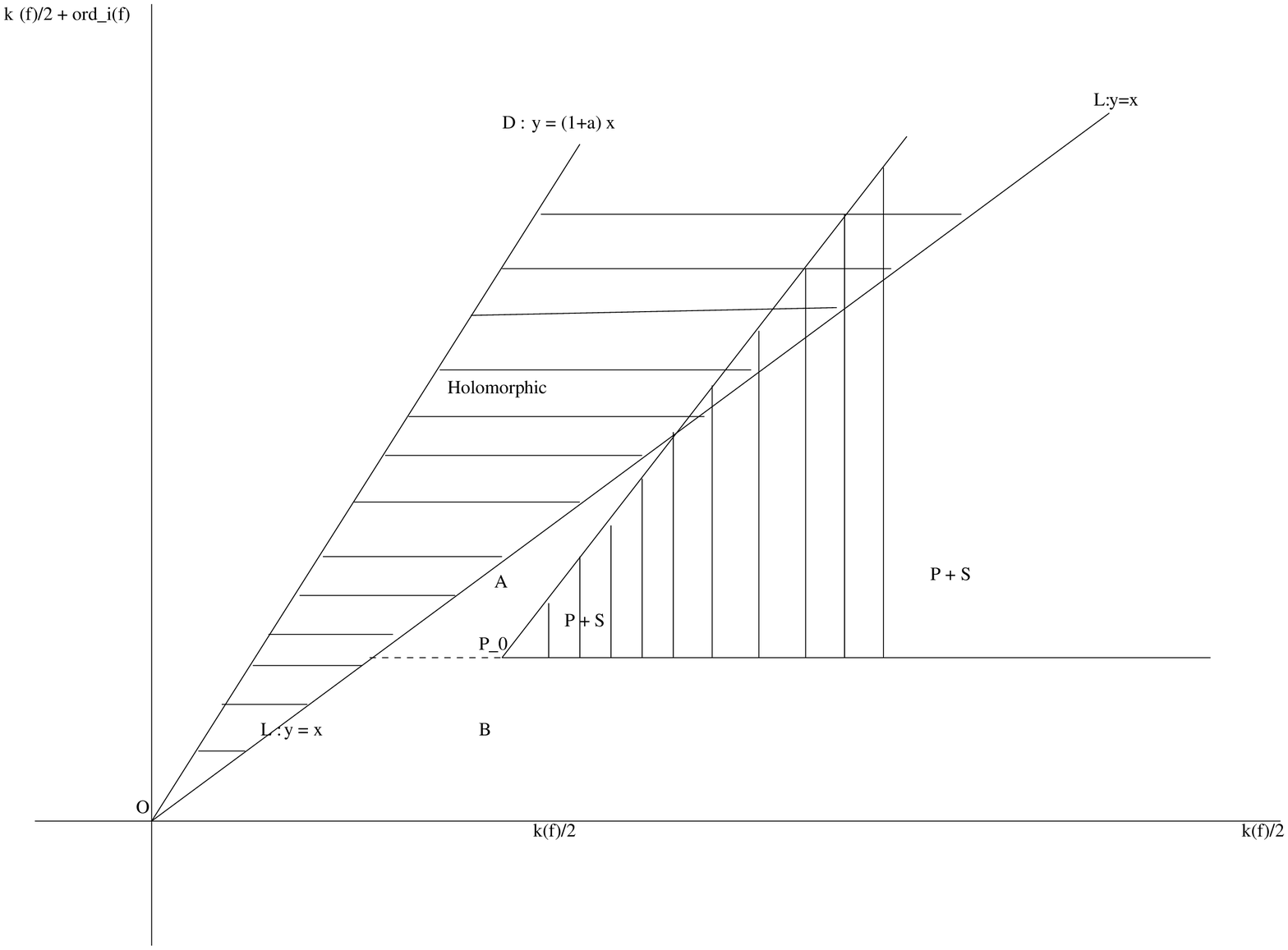} \vspace{2mm}
\end{figure}
here $k(f)$ is the weight of $f$  and $\rm{ord}_i(f)$ is
the  vanishing order of $f$ at $i.$ We define  $I(f)$ 
as the point $(I_1(f),I_2(f)),$ it is an invariant of $f.$ 

By Proposition~$6$, we have the property: 
$$ I(D_{\phi}(f)) = I(f) + (1,\beta),$$ 
with $\beta \geq 0.$ The case $\beta > 0$ may only happen 
if $I(f)$ is on the line $y=(2 \mathcal{K} +1) x.$

Let $J$ be the ideal of modular forms of 
positive weights over $\Gamma$ and $f_j,(j=1,
\cdots,d)$ be a basis of  $J / J^2.$
Let $R_0 = \langle f_1, \cdots,f_d, \omega , D_{\phi} (\omega) \rangle,$
the ring generated by $f_1,\cdots,f_d,\omega$ and $D_{\phi}(\omega).$ 
We will construct a sequence of subrings of $M_*(\Gamma; \{i \})$: 
$$ R_0 \subset R_1
\subset \cdots \subset R_i \subset R_{i-1} \subset \cdots $$ such that
for any $i$, $R_i$ is finitely generated. 
We will proof that this sequence is stationary
from a certain rank $n_0,$ and also that 
$D_{\phi}(R_{n_0}) = R_{n_0} = R_{n_0 +1}.$

We consider the semigroup of $\NM^2$ that is finitely  
generated by the set $I(R_0)=\{ I(f) | f \in 
R_0 \}$. For $f \in M_*(\Gamma),$ we have $\rm{ord}_i(f) 
\leq a \frac{k} {2}$ for a certain $a > 0$, as
a consequence of the   formula of zeros for modular forms.
If we choose $a$ to be  minimal, then the line $D$ defined by 
$y = (a + 1) x$ contains a non zero element of
$I(R_0)$ and $I(R_0)$ has $D$ as boundary. 
Hence, the semigroup $I(R_0)$ is embedded
in the sector $S$ delimitated  by the lines $(Ox)$ and $D.$ 
The intersection of  $I(R_0)$ and of the line $(Ox)$ contains 
$I(\omega) = (2,0)$ and 
$I(D_{\phi}(\omega)) = (3,0),$ so this intersection 
contains all points $(a,0)$ with $a \geq 3$.  
The set $I(R_0) \cap \{ (x,y) | y \geq x \}$ 
coincides with the set $I(M_*(\Gamma)).$ Indeed,
$\frac{k} {2} + \mathrm{ord}_i(f) \geq \frac{k} {2}$ 
if and only if $\rm{ord}_i(f) \geq 0.$ 
Furthermore the elements of $M_*(\Gamma ; \{ i \})$ 
have no poles outside the $\Gamma-$orbit of $i.$ 
It is clear that the group $I(R_0)$ is equal to $\ZM^2.$

We apply the lemma to the semigroup $I(R_0)$, 
and deduce from this that there exists $P_0 \in I(R_0)$
such that $ (P_0 + S) \cap \mathbb{Z}^2 \subset I(R_0),$
where $S$ is the sector associated to $I(R_0)$
like  in the Lemma. 

We have the essential property that: if $F$ is 
an element  of $M_*(\Gamma;\{ i \})$ and
$I(F) \in  (P_0 + S) \cap \ZM^2,$ then $F \in R_0$. 
Indeed there exists $g \in R_0$ such that
$I(F)=I(g).$ So there exists a linear combination
$g_1$ of $F$ and $g$ such that $I_1(g_1)=I_1(F)$
and $I_2(g_1) > I_2(F).$  
By reiterating this construction, we obtain
a sequence of  points $I(g_i)$ which for $i$ large 
exceed the line $y=x.$ 
By induction in the opposite direction, one deduces that $F \in R_0.$  
In particular, if $I(D_{\phi}(f_j)) \in (P_0 + S) \cap \ZM^2,$
then $D_{\phi}(f_j) \in R_0.$ 
  
There are two cases: either $I(f_j)$ is in the left of the sector $P_0 + S$
(region $A$ of the  diagram), and  then we must add the necessary number  
of derivations of $f_j$ to get in to the sector $P+S$.
Or case $f_j$  lies below the sector (region $B$
of the diagram): $I_2(f_j) < I_2(P_0)$. 

We define the set:
$$E(R_0) =\left\{ y \mid \nexists x \in \NM, (x,y) \in I(R_0) \right\}.$$ 
In other words, $E(R_0)$ is the set of horizontal lines 
not occupied by $I(R_0).$ 
We have $0 \not\in E(R_0),$ because $I(\omega) \in (Ox).$

There exists $y_0$ such that: if $x > y \geq y_0$ then $(x,y) \in I(R_0).$
There also exists $x_0$ such that if $ y < y_0$ and $y \not\in E(R_0)$ then
$I(R_0) \subset \{ (x,y) \mid x \geq x_0 \}.$ 

By induction, we define a sequence of rings:
$R_{j+1} = \langle R_j , 
D_{\phi}(R_j) \rangle.$ We also define a sequence of sets 
$E(R_0) \supset E(R_1) \supset \cdots,$  by:
$$E(R_j) =\left\{ y \mid \not\exists x \in \NM, (x,y) \in I(R_j) \right\}.$$
There exists $y_j$ such that if $x > y \geq y_j$ then $(x,y) \in I(R_j).$
There also  exists $x_j$ such that, if $ y < y_j$ and $y \not\in E(R_j)$ then
$I(R_j) \subset \{ (x,y) \mid x \geq x_j \}.$ Then, we have:
$$  \cdots \subset \cdots \subset E(R_{j+1}) \subset E(R_j) 
\subset \cdots \subset E(R_0).$$ The sequence of finite sets
$E(R_j)$ decreases, hence it is stationary. 
Then, there exists $j_0 \in \NM$ such that
$E(R_{j_0})= E(R_{j_0 + 1}).$ 
(this means, there will be no new  occupied lines). 

As a consequence, the sequence of rings $R_j$ is stationary 
from a  certain rank $n_0 \geq j_0$ with 
$D_{\phi}(R_{n_0}) = R_{n_0 }= R_{n_0 +1}.$ 
Indeed,  let $h_0 \in R_{j_0},$  after applying a  
suitable power of $D_{\phi}$ (this power is
bounded by certain integer $N$ independent of the
choice of $h_0$) we obtain $h_1 \in R_{n_0}$
 such that $I_1(h_1) \geq x_{n_0}$
(it is clear that $I_2(h_1) \not\in E(R_{n_0})$). Since there exists
$h_2 \in R_{n_0}$ such that $I(h_2)= I(D_{\phi}(h_1)),$ then
there exists a linear combination of $h_2$ and $D_{\phi}(h_1),$ 
equal to an element $h_3 \in R_{n_0 }$ such that its 
image under $I$ is a  point in the same vertical line as 
$I(D_{\phi}(h_1))$ and with strictly greater ordinate. 
By reiterating this construction, we obtain
a sequence of elements $(h_n)$ in $R_{n_0}.$ For $n$ enough large,
we have $I_2(h_n) \geq y_{n_0}$ so $h_n \in R_{n_0}.$ 
This proves that $D_{\phi}(R_{n_0}) = R_{n_0}$ and
$R_{n_0 + 1} = R_{n_0}.$ 

\end{proof}
\begin{remarque}
Let $\Gamma \subset \rm{PSL}(2, \RM)$ be a discrete  cocompact subgroup,
and let $R \subset M_*^{\mathrm{mer}}(\Gamma)$ be a  ring closed under
$D_{\phi}$ which  contains  $\omega.$  Then, the ring $\widetilde{R}=
R[\phi]$ is closed under $D.$
Indeed $D(f) = D_{\phi}f + k \phi f $ for $f\in R$ 
and $D(\phi)= \omega + \phi^2.$ 
\end{remarque}
\begin{corollaire}
{\bf Let $\Gamma \subset \rm{PSL}(2, \RM)$ be a discrete  cocompact 
subgroup. Then there exists a finitely generated ring
closed under $D$ which contains the ring of 
quasimodular forms over $\Gamma$.
One can take  $\widetilde{R} = R[\phi]$ with $\phi$ 
and $R$ such as in Theorem $10.$} 
\end{corollaire}
\begin{proof}
It is clear that $\widetilde{R}$ contains the ring of
modular forms $M_*(\Gamma)$. On the other hand, $\widetilde{R}$ 
is closed under $D$ by the last remark. By using Theorem $2,$
we deduce from this that $\widetilde{R}$ contains the ring 
$\widetilde{M}_*(\Gamma)$ of quasimodular forms. 
\end{proof}
\begin{remarque}
{\rm The ring $\widetilde{R}$ is generated by a finite number of 
meromorphic modular forms of positive
weights, so  in each weight
$k$, the $\mathbb{C}$ vectoriel space $\widetilde{R}_k$   
is finite dimensional.
Furthermore, we have $\dim \widetilde{R}_k  = \mathrm{O}(k^2),$ 
i.e. the finitely generated ring
$\widetilde{R}$ has the same order of growth
as its  non-finitely genertaed  subring $\widetilde{M}_*$.}
\end{remarque}
%%%%%%%%%%%%%%%%%%%%%%%%%%%%%%%%%%%%%%%%%%%%%%%%%%%%%%%%%%%%
%%%%%%%%%%%%%%%%%%%%%%%%%%%%%%%%%%%%%%%%%%%%%%%%%%%%%%%%%%%%%%%%%%%%%%%%%%%
%%%%%%%%%%%%%%%%%%%%%%%%%%%%%%%%%%%%%%%%%%%%%%%%%%%%%%%%%%%%%%%%%%%%%%%%%%%%

\end{document}